\documentclass{amsproc}

\newtheorem{thm}{Theorem}[section]
\newtheorem{cor}[thm]{Corollary}
\newtheorem{prop}[thm]{Proposition}

\theoremstyle{definition}
\newtheorem{defn}[thm]{Definition}

\title[M\"obius Transformations of the Circle Form a Maximal Conv. Group]
{M\"obius Transformations of the Circle Form a Maximal Convergence Group}
\author{Ara Basmajian}
\address{Department of Mathematics, University of Oklahoma, Norman, Oklahoma 73019}
\email{abasmajian@ou.edu}

\author{Mahmoud Zeinalian}
\address{Department of Mathematics, C.W.Post College of Long Island University, 720 Northern Boulevard, Brookville, New York 11548}
\email{mzeinalian@liu.edu}

\date{February 17, 2005 and, in revised form, August 15, 2005 }

\begin{document}
\bibliographystyle{h-elsevier2}

\maketitle

\begin{abstract}
We investigate the relationship between quasisymmetric and
convergence groups acting on the circle. We show that the M\"obius
transformations of the circle form a maximal convergence group.
This completes the characterization of the M\"obius group as a
maximal convergence group acting on the sphere. Previously,
Gehring and Martin had shown the maximality of the M\"obius group
on spheres of dimension greater than one. Maximality of the
isometry (conformal) group of the hyperbolic plane as a uniform
quasi-isometry group, uniformly quasiconformal group, and  as a
convergence group in which each element is  topologically
conjugate to an isometry may be viewed as consequences.
\end{abstract}

\section{Introduction}

The isometries of real hyperbolic space of dimension two or higher
induce conformal diffeomorphisms on its ideal boundary. In fact,
if the dimension is strictly greater than two, then all conformal
diffeomorphisms will arise in this way. In contrast, in dimension
two, conformality on the boundary is a trivial condition. For
instance, every diffeomorphism of a Riemannian circle preserves
the conformal class of the metric. It is for this reason that the
study of the group of M\"obius transformations of the circle
differs from its higher dimensional cousins. In the paper
\cite{G-M}, Gehring and Martin show that the M\"obius group is a
maximal convergence group acting on the boundary of real
hyperbolic space of dimension greater than two. This result was
extended (see \cite{BZ}) to the action of the isometry group of a
rank one symmetric space of noncompact type except the hyperbolic
plane.

In this note, we complete the characterization of the M\"obius
group as a maximal convergence group by considering  the remaining
case
 of the hyperbolic plane; namely, the group of M\"obius
transformations of the circle acts as a maximal convergence group
(Theorem \ref{thm:maxconv}). Other maximality statements, such as
the maximality of the isometry (conformal) group of the hyperbolic
plane as a uniform quasi-isometry group and a uniformly
quasiconformal group (see Corollary \ref{cor:psl(2,R)maximalqi}
and the discussion at the end of that section) may be regarded as
consequences of Theorem \ref{thm:maxconv}. See \cite{G-P} for
further discussion on quasi-isometry groups. Another implication
is the maximality of the isometry group of the hyperbolic plane as
a  convergence group in which each element is  topologically
conjugate to an isometry (Corollary
\ref{prop:psl(2,R)maximalconv}).

Let $X$ be a compact topological space. A family $\mathcal F$ of
orientation preserving homeomorphisms of $X$ is said to have the
{\it convergence property} if each infinite sequence $\{f_n\}$ of
$\mathcal F$ contains a subsequence which,
\begin{description}
\item[C1] converges uniformly to a homeomorphism of $X$, or
\item[C2]  has the attractor-repeller property, that is, there
exists a point $a \in X$, the {\it attractor}, and a point $r \in
X$, the {\it repeller}, so that the $\{f_n\}$ converge to the
constant function $a$, uniformly outside of any open neighborhood
of $r$. Note that $a$ may equal $r$.
\end{description}

We remark that the convergence groups considered in this paper are
 comprised only of orientation preserving homeomorphisms. We could
  equally as well include  orientation reversing
 homeomorphisms, in which case the theorems in this paper have obvious
 modifications that are left to the reader.
 Hence, whenever  homeomorphism is mentioned  in this paper it is assumed to be
 orientation preserving.

\section{Elementary facts about  quasisymmetric mappings}

In this section, we assemble  some  elementary  facts which will
be needed later in the paper. For the basics on quasisymmetric and
quasiconformal maps, we refer to the following papers and books:
\cite{A}, \cite{D-E}, \cite{G-L}, \cite{H}, \cite{L}, and
\cite{V}. For M\"obius groups and  hyperbolic geometry, the reader
may consult \cite{Be} or \cite{M}.

Let $\mathbb{H}$ denote the upper half plane, endowed with the
hyperbolic metric, and $\widehat{\mathbb{R}}=\mathbb{R} \cup
\{\infty\}$ denote its ideal boundary. $\widehat{\mathbb{R}}$ can
be identified with the unit circle $S^1$. The group $\text{PSL}(2,
\mathbb{R})=\{z\mapsto \frac{az+b}{cz+d}: a, b, c, d \in
\mathbb{R},  ad-bc=1\}$ is the full group of orientation
preserving isometries of $\mathbb{H}$. This group is also the full
group of conformal homeomorphisms of $\mathbb{H}$. Let
$\text{M\"{o}b}^+(\widehat{\mathbb{R}})$ denote the group of
homeomorphisms of $\widehat{\mathbb{R}}$ which are induced by the
isometries of $\mathbb{H}$. Note that $\text{PSL}(2, \mathbb{R})$
and $\text{M\"{o}b}^+(\widehat{\mathbb{R}})$ are  isomorphic
groups which act on different spaces.

When dealing with mappings of the circle we will need to normalize
by post composition using an element of
$\text{M\"{o}b}^+(\widehat{\mathbb{R}})$   so that infinity is a
fixed point of the map. Observe that an element of the stabilizer
of $\infty$ in $\text{M\"{o}b}^+(\widehat{\mathbb{R}})$ is a
linear map, $x\mapsto mx+b$, where $m>0$ and $b \in \mathbb{R}$.
Hence given a symmetric configuration of point triples,
$\{x-t,x,x+t\}$, its image remains a symmetric configuration.

\begin{defn}\label{defn:quasisymmetric}   Let
$f:\widehat{\mathbb{R}}\rightarrow \widehat{\mathbb{R}}$ be an
orientation preserving homeomorphism and $k>0$. The homeomorphism
$f$ is called $k$-quasisymmetric if after normalizing so that it
fixes infinity, $f$ satisfies,
$$\frac{1}{k} \leq \frac{f(x+t)-f(x)}{f(x)-f(x-t)}\leq
k$$ for all $x \in \mathbb{R}$ and all $t>0$. In other words, the
image of equal length juxtaposed intervals has uniform bounded
length ratio.
\end{defn}

Given the observation preceding  Definition
\ref{defn:quasisymmetric}, it is easy to see that the  condition
of being $k$-quasisymmetric is independent of the normalizing
M\"obius transformation, the elements of
$\text{M\"{o}b}^+(\widehat{\mathbb{R}})$ are $1$-quasisymmetric,
and that post or precomposition by linear maps does not change the
quasisymmetric constant.

 In the sequel, we will need the fact that
$\text{M\"{o}b}^+(\widehat{\mathbb{R}})$ is the full group of
$1$-quasisymmetric homeomorphisms. To see this, using the triple
transitivity of $\text{M\"{o}b}^+(\widehat{\mathbb{R}})$, it is
enough to show that  a   $1$-quasisymmetric mapping $f$ which
fixes $0,1,$ and $\infty$ is the identity. Now, the triple $\{0,
\frac{1}{2}, 1\}$ must be taken to a symmetric triple, hence $f$
also fixes $\frac{1}{2}$. Similar considerations  allow us to
conclude that all  rationals of the form $\frac{n}{2^k}$   are
fixed. Finally, using continuity, $f$ must fix every real number
and thus is the identity homeomorphism.

The following proposition is a classical result. See Ahlfors
\cite{A} or Lehto \cite{L} for reference.

\begin{prop}\label{extension2} Let
$f:\widehat{\mathbb{R}}\rightarrow \widehat{\mathbb{R}}$ be an
increasing  homeomorphism with $f(\infty)=\infty$. Assume that
$$\frac{1}{k} \leq \frac{f(x+t)-f(x)}{f(x)-f(x-t)}\leq
k$$ for all $x$ and all $t>0$. Then there exists a K(k)
-quasiconformal homeomorphism of $\mathbb{H}$ which extends $f$.
The number $K(k)$ depends only on $k$.
\end{prop}

\section{Maximality of $\text{M\"{o}b}^+(\widehat{\mathbb{R}})$}

A family $\mathcal{F}$ is said to be {\it uniformly quasisymmetric
(quasiconformal)} if all the maps in $\mathcal{F}$ are
$k$-quasisymmetric ($K$-quasiconformal) for some $k$ (for some
$K$).

\begin{prop} \label{prop:conv} A uniformly quasisymmetric family $\mathcal{F}$ of
homeomorphisms of $\widehat{\mathbb{R}}$ is a convergence family.
In particular, $\text{M\"{o}b}^+(\widehat{\mathbb{R}})$ is a
convergence group.
\end{prop}

\begin{proof} Using Proposition \ref{extension2}, there exists  a number $K$ such that every element of this family can be extended to a $K$-quasiconformal mapping of $\mathbb{H}$.
Using the fact that a sequence of distinct $K$-quasiconformal
mappings of $\mathbb{H}$ has a subsequence which either converges
to a $K$-quasiconformal map or has the attractor-repeller property
with attractor and repeller on the boundary (see \cite{V},
Corollaries 19.3 and 37.4, or extend each map in $\mathcal{F}$ to
$\text{S}^2$ by reflection and use the results of \cite{G-M}), we
may conclude it acts as a convergence family on $S^1$.
\end{proof}

\begin{prop}\label{prop:uniformlyqs} Let $\mathcal{F}$ be a family of homeomorphisms of
$\widehat{\mathbb{R}}$ which is closed under post and
precomposition by elements of
$\text{M\"{o}b}^+(\widehat{\mathbb{R}})$. Then $\mathcal{F}$ has
the convergence property if and only if $\mathcal{F}$ is a
uniformly quasisymmetric family.
\end{prop}

\begin{proof} Suppose $\mathcal{F}$ has the convergence property.
Let \begin{equation} \mathcal{F}^{\prime}=\{f\in
\mathcal{F}:f(0)=0,f(1)=1,\,\text{and}\,f(\infty)=\infty \}.
\end{equation} Since $\mathcal{F}^{\prime}$ has the convergence property, it must
be that there are negative constants $M$ and $m$ so that,
$m<f(-1)<M<0$, for all $f \in \mathcal{F}^{\prime}$. Any element
of $\mathcal{F}$ can be post composed by an element of
$\text{M\"{o}b}^+(\widehat{\mathbb{R}})$ to yield an element of
$\mathcal{F}^{\prime}$. Since any triple $\{x-t, x, x+t\}$ in
$\mathbb{R}$ can be moved by Euclidean translation and dilation to
$\{-1,0,1\}$, we may conclude that the elements of $\mathcal{F}$
form a uniformly quasisymmetric family. The converse follows from
Proposition \ref{prop:conv}.

\end{proof}

\begin{thm} \label{thm:maxconv} $\text{M\"{o}b}^+(\widehat{\mathbb{R}})$ acts on
$\widehat{\mathbb{R}}$ as a maximal convergence group. That is,
there is no convergence group that properly contains
$\text{M\"{o}b}^+(\widehat{\mathbb{R}})$.
\end{thm}

\begin{proof} The fact that
$\text{M\"{o}b}^+(\widehat{\mathbb{R}})$ is a convergence group
follows from Proposition \ref{prop:conv}. Next let $G$ be a
convergence group acting on $\widehat{\mathbb{R}}$ containing
$\text{M\"{o}b}^+(\widehat{\mathbb{R}})$. Using Proposition
\ref{prop:uniformlyqs}, we know that the action of $G$ is as a
uniformly quasisymmetric group. On the other hand, suppose there
exists an element $g \in G$ not contained in
$\text{M\"{o}b}^+(\widehat{\mathbb{R}})$. This means that after
normalizing $g$, so that it fixes $\infty$, there must be three
symmetrically spaced points in $\mathbb{R}$ where the
quasisymmetric constant is not $1$. Post and precomposing by
Euclidean translations and dilations, we may assume that the three
points are $\{-1,0,1\}$ and that  $g$ fixes $0$ and $1$. Since
Euclidean translation and dilation do not effect the
quasisymmetric constant for a triple, it must be that $g$ does not
fix $-1$. By possibly replacing $g$ with $g^{-1}$, we may assume
that $g$ takes $-1$ into the interval $(-1,0)$. Clearly $g^n(-1)$
is an increasing sequence of negative numbers and hence has a
limit $y$ which is necessarily a fixed point of $g$. Since  $<g>$
is  a convergence group, $y$ is strictly less than $0$. Next
consider the triple of points $\{-1,y,0\}$. The length ratio of
the juxtaposed intervals $[-1,y]$ and $[y,0]$ is,

\begin{equation}
\frac{|0-y|}{|y-(-1)|}=\frac{|y|}{|y+1|}.
\end{equation}
The length ratio of image intervals under the iterates of $g$ are
\begin{equation}
\frac{|g^n(0)-g^n(y)|}{|g^n(y)-g^n(-1)|}=\frac{|y|}{|y-g^n(-1)|}
\end{equation} which goes to $\infty$, as $n\rightarrow \infty$. This contradicts
the fact that $G$ is uniformly quasisymmetric. Hence, it must be
that the quasisymmetric  constant for $g$ is $1$, and thus $g \in
\text{M\"{o}b}^+(\widehat{\mathbb{R}})$.
\end{proof}

\section{Maximality of $\text{PSL}(2, \mathbb{R})$}

An immediate corollary of Theorem \ref{thm:maxconv} is,

\begin{cor}
The M\"obius group, $\text{M\"{o}b}^+(\widehat{\mathbb{R}})$, is a
maximal uniformly quasisymmetric group.

\end{cor}

 A homeomorphism $f: \mathbb H \rightarrow \mathbb H $ is said to
be a {\it quasi-isometry} if there exist positive constants $A$
and $B$ so that \begin{equation} A^{-1} d (x_1, x_2) -B \leq d(f
(x_1), f(x_2)) \leq A d (x_1, x_2) +B
\end{equation}
for all $x_1,x_2 \in \mathbb{H}$. The constant $A$ is referred to
as the {\it Lipschitz constant} of the quasi-isometry. For a
general reference on quasi-isometries, we refer the reader to
\cite{G-P}. It is well known that a quasi-isometry continuously
extends to the boundary and that the induced map on the boundary
is a quasisymmetric homeomorphism. One defines an equivalence
relation on quasi-isometries by declaring two  to be equivalent if
they induce the same homeomorphism on the boundary. Let
$\text{QI}(\mathbb{H})$ denote the group of equivalence classes of
quasi-isometries of $\mathbb{H}$. Since the natural map from
$\text{PSL}(2, \mathbb{R})$ into $\text{QI}(\mathbb{H})$ is
injective, we will continue to denote its image with the same
notation.  A family ${\mathcal{F}} \subset \text{QI}(\mathbb{H})$
is said to be a {\it uniformly quasi-isometric} family if each
equivalence class has  Lipschitz constant  less than a uniform
bound. The following is a simple consequence of Theorem
\ref{thm:maxconv}, observed in Gromov and Pansu (see \cite{G-P}).

\begin{cor}\label{cor:psl(2,R)maximalqi}
 Let $G \leq \text{QI}(\mathbb{H})$ be a uniform quasi-isometry group acting on
the hyperbolic plane $\mathbb{H}$. If $\text{PSL}(2,
\mathbb{R})\leq G$, then $G=\text{PSL}(2, \mathbb{R})$.
\end{cor}

\begin{proof} A quasi-isometry $f$ of $\mathbb{H}$ extends to a homeomorphism of
$\mathbb{H} \cup \widehat{\mathbb{R}}$. Moreover, the induced
mapping $f|_{\widehat{\mathbb{R}}}$ on $\widehat{\mathbb{R}}$ is a
quasisymmetric mapping where the quasisymmetric constant  depends
only on the Lipschitz constant of the quasi-isometry. Consider the
homomorphism $\phi: G \rightarrow Homeo(\widehat{\mathbb{R}})$,
given by $[f] \mapsto f|_{\widehat{\mathbb{R}}}$. Note that
$\text{Image}(\phi)$ is a uniformly quasisymmetric group which
contains $\text{M\"{o}b}^+(\widehat{\mathbb{R}})$. By Proposition
\ref{prop:conv}, $\text{Image}(\phi)$ is a convergence group.
Since $\text{M\"{o}b}^+(\widehat{\mathbb{R}})$ is a maximal
convergence group (Theorem \ref{thm:maxconv}),
$\text{Image}(\phi)$ equals
$\text{M\"{o}b}^+(\widehat{\mathbb{R}})$. Injectivity of $\phi$ is
a tautology. Since $\text{Image}(\phi)=\phi(\text{PSL}(2,
\mathbb{R}))= \text{M\"{o}b}^+(\widehat{\mathbb{R}})$, we conclude
that $G=\text{PSL}(2, \mathbb{R})$.
\end{proof}

As in the  proof above, it is easy to see that  a convergence
group acting on $\mathbb{H}\cup\widehat{\mathbb{R}}$ which
contains $\text{PSL}(2, \mathbb{R})$ induces the action of
$\text{M\"{o}b}^+(\widehat{\mathbb{R}})$ on the boundary,
$\widehat{\mathbb{R}}$. Furthermore, if each element of this
convergence group is topologically conjugate to an element of
$\text{PSL}(2, \mathbb{R})$, then the induced action has trivial
kernel. Hence the natural homomorphism given by restriction of the
convergence group to the boundary is in fact an isomorphism onto
$\text{M\"{o}b}^+(\widehat{\mathbb{R}})$. Since the image of
$\text{PSL}(2, \mathbb{R})$ is
$\text{M\"{o}b}^+(\widehat{\mathbb{R}})$, it must be that the
convergence group equals $\text{PSL}(2, \mathbb{R})$. We have
proven,

\begin{cor}\label{prop:psl(2,R)maximalconv}
 Let $G$ be a convergence group acting on $\mathbb{H}\cup
\widehat{\mathbb{R}}$. Suppose that every element of $G$ is
topologically conjugate to an element of $\text{PSL}(2,
\mathbb{R})$. If $\text{PSL}(2, \mathbb{R})\leq G$, then
$G=\text{PSL}(2, \mathbb{R})$. The conjugating homeomorphism need
not be the same for all elements of $G$.
\end{cor}

The reader should compare the above corollary to the fact that
$\text{PSL}(2, \mathbb{R})$ is a maximal uniformly quasiconformal
group acting on $\mathbb{H}$. That is, if $G$ is a uniformly
quasiconformal group containing $\text{PSL}(2, \mathbb{R})$, then
$G=\text{PSL}(2, \mathbb{R})$. This fact follows from the results
of \cite{S} and \cite{T}.

\bibliographystyle{amsalpha}

\begin{thebibliography}{A}
\bibitem[A]{A} L. V. Ahlfors, \textit{Lectures on Quasiconformal Mappings}, Van Nostrand-Reinhold, Princeton, New Jersey,
(1966).

\bibitem[Be]{Be}Alan F. Beardon, \textit{The geometry of discrete groups},
 Graduate Texts in Mathematics, 91.
Springer-Verlag, New York, 1983. xii+337 pp.

\bibitem[B-Z]{BZ} A. Basmajian, M. Zeinalian, \textit{Maximal
convergence groups and rank one symmetric spaces}, preprint.

\bibitem[D-E]{D-E} A. Douady, C. J. Earle, \textit{Conformally natural extension of homeomorphisms of the
circle}, Acta Math. 157 no. 1-2 (1986), 23--48.

\bibitem[Ga-L]{G-L} F. P. Gardiner, N. Lakic, \textit{Quasiconformal
Teichmuller Theory}, American Mathematical Society, Mathematical
Surveys and Monographs, 76 (2000).

\bibitem[Ge-M]{G-M} F. W. Gehring, G. J. Martin, \textit{Discrete quasiconformal groups {I} },
Proc. London Math. Soc. (3), Vol 55, (1987) 331--358.




\bibitem[Gr-P]{G-P} M. Gromov, P. Pansu, \textit{Rigidity of lattices: An
introduction}, Geometric Topology: Recent developments, Lecture
Notes in Math. 1504, Springer, Berlin (1991), 39--137.

\bibitem[H]{H}  A. Hinkkanen, \textit{Uniformly quasisymmetric groups}.
Proc. London Math. Soc. (3) 51 (1985), no. 2, 318--338.

\bibitem[L]{L} O. Lehto, \textit{Univalent Functions and Teichm\"{u}lar Spaces}, Springer-Verlag, (1987).



\bibitem[M]{M}B. Maskit, \textit{Kleinian groups}, Grundlehren der Mathematischen
Wissenschaften [Fundamental Principles of Mathematical Sciences],
287. Springer-Verlag, Berlin, 1988. xiv+326 pp.



\bibitem[S]{S} D. Sullivan, \textit{On the ergodic theory at infinity of an arbitrary discrete
group of hyperbolic motions}, Riemann surfaces and related topics:
Proceedings of the 1978 Stony Brook Conference (1978), 465--496;
Ann. of Math. Stud., 97, Princeton Univ. Press (1981).

\bibitem[T]{T} P. Tukia, \textit{On two-dimensional quasiconformal groups},
 Ann. Acad. Sci. Fenn. Ser. A I Math. 5 no. 1 (1980), 73--78.

\bibitem[V]{V} J. Vaisala, \textit{Lectures on $n$-dimensional quasiconformal
mappings}, Lecture Notes in Math. 229, Springer-Verlag, Berlin-New
York (1971).
\end{thebibliography}

\author{Ara Basmajian, Department of Mathematics, University of Oklahoma, Norman, OK
73019, abasmajian@ou.edu}

\author{Mahmoud Zeinalian, Department of Mathematics, C.W.Post Campus, Long Island University, 720 Northern Boulevard, Brookville, NY 11548, mzeinalian@liu.edu}
\end{document}